\documentclass[12pt]{article}
\usepackage{amsmath}
\usepackage{amssymb}
\usepackage{amscd}
\usepackage{epsfig}

\def\no{\noindent}
\numberwithin{equation}{subsection}

\newenvironment{theorem}{\par\bigskip\noindent\refstepcounter{equation}{\bf Theorem \theequation.}\em}{\em\par\bigskip\noindent}

\newcommand{\Z}{\mathbb Z}

\newcommand{\qed}{{\unskip\nobreak\hfil
        \penalty50\hskip1em\hbox{}\nobreak\hfil
        $\square$\parfillskip=0pt\finalhyphendemerits=0 \par}}
\newcommand{\proof}{\no{\em Proof.\ }}

\def\ga{\gamma}

\begin{document}

\title{The volume and lengths on a three sphere}
\author{Christopher B. Croke\thanks{Supported by NSF grant DMS-99-71749.}}
\date{\today}
\maketitle

\begin{abstract}
We show that the volume of any Riemannian metric on a three sphere is bounded below by the length of the shortest closed curve that links its antipodal image.  In particular, the volume is bounded below by the minimum of the length of the shortest closed geodesic and the minimal distance between antipodal points.
\end{abstract}

\setcounter{section}{1}
\setcounter{subsection}{0}

\subsection{Introduction}
In this paper we consider the volume, $Vol(g)$, of a Riemannian metric $g$ on the three sphere. We let $L(g)$ represent the length of the shortest nontrivial closed geodesic in $(S^3,g)$.  By an ``antipodal map'', $A$, on an $n$-sphere we will mean an order 2, fixed point free diffeomorphism.  We will let $D(g,A)= \inf_{x\in S^3} d(x,Ax)$, and $D(g) = sup_A D(g,A)$.  The main result of this paper is:

\begin{theorem}
\label{main}
For any Riemannian metric on $S^3$ we have:
$$Vol(g)^{\frac 1 3}\geq C_1 \min\{L(g),2D(g)\}$$
where $C_1$ is a universal constant (which can be taken to be $\frac 1 {3180}$)
\end{theorem}

This theorem is an example of a ``universal'' inequality.  This is a term introduced by Berger to describe inequalities on a Riemannian manifold $(M,g)$ (usually between minimizing objects in some topological class) which hold for all metrics $g$ on $M$.  With the glaring exception of Gromov's work \cite{Gr} (which we discuss below) most such inequalities are known in 2 dimensions only.  Many of these inequalities involve the systole, $sys(g)$, which is the length of the shortest closed noncontractible curve in $M$.  Estimates of the form $A(g)^{\frac 1 2}\geq c(M)sys(g)$, Where $A(g)$ represents the area, have been proved for all surfaces $M$.  The first of these was proved by Loewner (unpublished) where he finds the sharp value of $c(T^2)$ (sharp for the flat regular hexagonal torus).  The sharp values of $c(RP^2)$ and $c(K^2)$ for the projective plane (the round one is best) and Klein bottle (the best one is singular and not flat!) were proven by Pu \cite{Pu} and Bavard \cite{Ba} respectively.  For surfaces of genus $\gamma$, nonsharp constants $c(\gamma)$ (which unfortunately went to 0 as $\gamma$ grew) were first proved independently by Accola \cite{Ac} and Blatter \cite{Bl}.  Gromov greatly improve these constants in \cite{Gr} so that they grow correctly with $\gamma$.  However, the best constants are still unknown.  Finding them is likely to be an extremely hard problem.

None of the above says anything about $S^2$ since it is simply connected.  However there are two inequalities of this type on $S^2$.  The first, due to Berger \cite{Be} says that there is a universal constant $c_2$ such that for any antipodal map $A:S^2\to S^2$:
$$Vol(S^2,g)^{\frac 1 2}\geq c_2 D(g,A).\leqno{1.1}$$ 
The sharp constant is not known, but is conjectured to be that of the round sphere.  The other result (see \cite{Cr2}) is that 
$$Vol(S^2,g)^{\frac 1 2}\geq c_3 L(g).\leqno{1.2}$$
Also in this case the sharp constant is not known.  It is conjectured to be that for the singular metric one gets by gluing two flat equilateral triangles along their boundaries. 

We now turn to higher dimensions.  The most important results in the area are in Gromov's paper ``Filling Riemannian Manifolds'' \cite{Gr}.  Gromov shows that for essential manifolds (which includes $T^n$, $RP^n$, and all compact $K(\pi,1)$ spaces) we have $V(g)^{\frac 1 n}\geq c(n)sys(g)$.  However, this says nothing directly about $S^n$ (or any compact simply connected manifold).  For example the natural generalizations of 1.1 and 1.2 above for $S^2$ are open questions for $S^n$, $n\geq 3$.  Note that for any Riemannian manifold $(M,g)$ of injectivity radius $inj(M)$ and any $r\leq \frac 1 2 inj(M)$ the metric spheres $S(x,r)$ with their induced Riemannian metrics $\bar g$ satisfy $D(\bar g,A)\geq 2r$, so an estimate $Vol(g)^{\frac 1 n}\geq c(S^n)D(g,A)$ which was sharp for the round spheres would give sharp estimates for the volume of small metric spheres.  See \cite{Cr1} for nonsharp estimates for the volume of such metric spheres. 

In the other direction, Ivanov (see \cite{I1} or \cite{I2}) has given examples of a sequence of metrics on $S^3$ that Gromov-Housdorff converge to the standard metric but whose volumes go to zero. 

Although Theorem \ref{main} does not yield either 1.1. or 1.2 for $S^3$ it does show that for any given metric one or the other inequality must hold. 

For $A:S^3\to S^3$ an antipodal map we let $L(g,A)$ be the infimum of the lengths of curves, $\ga$ such that $\ga$ links $A(\ga)$.  Here we will say $\ga$ links $\tau$ if either $\ga$ intersects $\tau$ or the linking number of $\ga$ and $\tau$ is nonzero.  We consider two cases.  If there is no closed curve $\ga$ that links $A(\ga)$ of length less than $2D(g,A)$ then $L(g,A)=2D(g,A)$ and the minimum length $\ga$ is the union of two minimizing geodesics between $x$ and $A(x)$ for some $x$. Otherwise there is such a closed curve $\ga$.  In this case we can apply curve shortening to $\ga$ to get a continuous family $\ga_t$ each of which is shorter than $\ga$ and hence does not intersect $A(\ga_t)$.  Thus we see that $\ga_t$ links with $A(\ga_t)$ and hence has length bounded below by $L(g,A)$.  By looking at limits of such $\ga_t$ we find nontrivial closed geodesics of length less than the length of $\ga$.  Since this is true for all such $\ga$ we find that there is a closed geodesic of length $L(g,A)$.  Thus the invariants $L(g)$ and $D(g,A)$ are thus related to $L(g,A)$ via:
$$L(g,A)\geq min\{L(g),2D(g,A)\}.$$

Theorem \ref{main} thus will follow from:

\begin{theorem}
\label{main2}
For any Riemannian metric on $S^3$ and any antipodal map $A$, we have:
$$Vol(g)^{\frac 1 3}\geq C_1 L(A,g)$$
where $C_1$ is a universal constant (which can be taken to be $\frac 1 {3180}$)
\end{theorem}

The fundamental result in the proof of Gromov's isosystolic inequality is his filling radius theorem: $Fillrad(g)\leq c_n Vol(g)^{\frac 1 n}$, which holds for all Riemannian $n$-manifolds (here $c_3<265$).  We define the filling radius, $Fillrad(g)$, in the next section.  We also make fundamental use of it since Theorem \ref{main2} will thus follow from

\begin{theorem}
\label{fill3}
For any Riemannian metric on $S^3$ and any antipodal map $A$, we have:
$$Fillrad(g)\geq C_4 L(g,A)$$
Where $C_4$ is a universal constant (which can be taken to be $\frac 1 {12}$) 
\end{theorem}

The proof of this theorem is a generalization of an argument that works on $S^2$ to yield:

\begin{theorem}
\label{fill2}
For any Riemannian metric on $S^2$ we and any antipodal map $A$ we have:
$$Fillrad(g)\geq \frac 1 4 D(g,A)$$
\end{theorem}

This theorem along with the Filling radius theorem recovers Berger's estimate equation 1.1 (although with a worse constant).  We present the argument in section \ref{secfill2}. 

The author would like to thank Herman Gluck for helpful conversations.

\subsection{Notation/Preliminaries}

The main purpose of this section is to set up notation and remind the reader about the Filling radius (for details see \cite{Gr}):

In \cite{Gr} Gromov introduced the notion of the filling radius (we use integer coefficients) $Fillrad(M)$ of a closed $n$-dimensional manifold $M$ with a metric $d$ (not necessarily Riemannian).  We will only consider the case where $M$ is homeomorphic to $S^2$ or $S^3$ and the metric is Riemannian.  There is a natural isometric embedding (in the metric space sense!) $i:M\to L^\infty(M)$ defined by $i(x)(\cdot)=d_M(x,\cdot)$.  The filling radius is the infimum of $r$ such that $i(M)$ bounds in the tubular neighborhood $T_r(i(M))$ in the sense that $i_*(H_n(M;\Z))$ vanishes in $H_n(T_r(i(M));\Z)$.  We will represent a filling as continuous map $\sigma:\Sigma \to T_r(i(M))$ from an $n+1$-dimensional simplicial complex $\Sigma$ such that $\sigma|_{\partial \Sigma}:\partial \Sigma\to i(M)$ represents a generator in $H_n(M;\Z)$.  

We note that for any fixed $\epsilon>0$ by taking Barycentric subdivisions as needed we may assume that the $\sigma$-image of any simplex has diameter less than $\epsilon$ in $L^\infty(X)$.

We note that there can be no continuous map $f:\Sigma\to i(X)$ which agrees with $\sigma$ on $\partial\Sigma$ since $\sigma|_{\partial\Sigma}$ represents a generator of the top homology (so is not a boundary).  Our proof of Theorem \ref{fill2} will be by contradiction.  We assume that $Fillrad(g)$ is small, take a filling $\sigma:\Sigma\to L^\infty(M)$ as above, and show that $\sigma:\partial\Sigma\to i(M)$ extends to a continuous $f:\Sigma\to i(M)$ giving the desired contradiction.  Our proof of Theorem \ref{fill3} uses a similar contradiction.

We use the notation $S_i(\Sigma)$ to denote the $i$-skeleton of a simplicial complex $\Sigma$.  For each $i=0,1,2,...,n+1$  We will let $\{\Delta^i_j\}$ denote the set of $i$ simplices in the $n+1$ dimensional simplicial complex $\Sigma$.  For each i-simplex $\Delta^i_j$ we let $G^i_j$ be the connected graph (i.e. 1-complex) $G^i_j=\cup \{S_1(\Delta^{n+1}_k)| \Delta^i_j$ is a face of $\Delta^{n+1}_k\}$.  By taking Barycentric subdivisions if needed we can assume that any two simplices intersect in a single (possibly empty) common face. Then we see that $H_1(G^i_j)$ is generated by the boundaries of the $2$-simplices that are faces of the $n+1$ simplices in the above union, since Van Kampen's theorem implies that the two complex $\cup \{S_2(\Delta^{n+1}_k)| \Delta^i_j$ is a face of $\Delta^{n+1}_k\}$ is simply connected.

\subsection {Proof of Theorem \ref{fill2}}
\label{secfill2}

\proof (of Theorem \ref{fill2})
Let $\epsilon>0$ and choose a filling $\sigma:\Sigma\to L^\infty(S^2,g)$ of $(S^2,g)$ in the $Fillrad(g)+\epsilon$ tubular neighborhood of $S^2\subset L^\infty(S^2,g)$ (we will confuse $S^2$ with $i(S^2)$ since $i$ is an isometric embedding).  By taking subdivisions we can assume that $diam(\sigma(\Delta^3_k))<\epsilon$ for each $\Delta^3_k\in \Sigma$.  We will prove the theorem as suggested in the previous section by 
showing that if $Fillrad(g)<\frac {D(g,A)} {4}$ then we can find a continuous map $f:\Sigma \to S^2$ extending the map $\sigma|_{\partial \Sigma}$.  We choose $\epsilon$ so small that $2Fillrad(g)+3\epsilon\leq \frac{D(g,A)}{2}$

{\bf Step1: The 0 and 1 Skeletons}

We define $f$ on the 0-skeleton, $S_0(\Sigma)$, of $\Sigma$ by mapping each 0 simplex $v$ to a point $f(v)$ on $i(S^2)$ closest to $\sigma(v)$, hence
$$d_{L^\infty(S^2)}(f(v),\sigma(v)) < Fillrad(g)+\epsilon.$$
In particular, $f$ takes 0-simplices in the boundary to the same point as $\sigma$ does.  We note that if two vertices $v_1$ and $v_2$ are the endpoints of an edge then $$d_{S^2}(f(v_1),f(v_2))=d_{L^\infty(S^2)}(f(v_1),f(v_2)) \leq $$
$$\leq d_{L^\infty(S^2)}(f(v_1),\sigma(v_1))+d_{L^\infty(S^2)}(\sigma(v_1),\sigma(v_2))+d_{L^\infty(S^2)}(\sigma(v_2),f(v_2))< $$
$$< 2Fillrad(g)+3\epsilon\leq \frac{D(g,A)}{2}.$$
We define $f$ to map each nonboundary $\Delta^1_i\in S_1(\Sigma)$  to a minimizing geodesic between the $f$ image of the endpoints, hence the length $L(\Delta^1_i)$ satisfies 
$L(\Delta^1_i)<  \frac{D(g,A)}{2}$.
For $\Delta^1_i$ on the boundary, we let $f|_{\Delta^1_i}=\sigma|_{\Delta^1_i}$ and hence $Diam(f(\Delta^1_i))<\epsilon$.

{\bf Step 2: The 2 skeleton}

We now extend $f$ to $S_2(\Sigma)$.   We note that $f(\partial\Delta^2_j)\cap A\circ f(G^2_j)=\emptyset$, for otherwise there would be a point $x\in f(G^2_j)$ such that $A(x)\in f(G^2_j)$ but this cannot happen since step 1 guarantees that the diameter of $G^2_j$ is less than $D(g,A)$. 

We can thus extend $f$ to $S_2(\Sigma)$ in such a way that we have 
$$f:\Delta^2_j\to S^2 - A\circ f(G^2_j).$$
Note that for boundary simplices, $\Delta^2_j$, this will hold when we take $f|_{\Delta^2_j}=\sigma|_{\Delta^2_j}$ by the triangle inequality.  

{\bf Step 3: The 3 skeleton}

Now for every $\Delta^3_k$ in $\Sigma$ the previous step guarantees that $$f:\partial \Delta^3_k\to S^2-A\circ f(S_1(\Delta^3_k)).$$  
Hence, since $f|_{\partial \Delta^3_k}$ misses a point, we can extend $f$ to $\Delta^3_k$. 

This completes the proof.

\subsection {Proof of Theorem \ref{fill3}}

Throughout this section $S^3$ will be endowed with a fixed Riemannian metric $g$ and an antipodal map $A:S^3\to S^3$.

\bigskip
We will proceed analogously to the proof of Theorem \ref{fill2}.  However, in this case we will get our contradiction by finding a singular chain in $S^3$ whose boundary is $\sigma$ restricted to $\partial \Sigma$.  We will do this one skeleton at a time, associating to each simplex $\Delta_i^j$ a singular simplicial chain $c_i^j$ whose boundary $\partial c_i^j$ corresponds to the already defined chain associated to $\partial \Delta_i^j$ (i.e. if $\partial \Delta_i^j = \Sigma_k (-1)^{\alpha(i,k)}\Delta_k^{j-1}$ then $\partial c_i^j = \Sigma_k (-1)^{\alpha(i,k)}c_k^{j-1}$).
We do this while associating simplices $\Delta^j_i$ of $\partial \Sigma$ to the chain consisting only of $\sigma$ applied to $\Delta^j_i$.

Let $\epsilon>0$ and choose a filling $\sigma$ of $(S^3,g)$ in the $Fillrad(g)+\epsilon$ tubular neighborhood of $S^3\subset L^\infty(S^3,g)$.  By taking subdivisions we can assume that $diam(F(\Delta^4_l))<\epsilon$ for each $\Delta^4_l\in \Sigma$.  We will prove the theorem by 
showing that if $Fillrad(g)<\frac {L(g,A)} {12}$ then we find a chain as above extending $\sigma|_{\partial \Sigma}$. We choose $\epsilon$ so small that $2Fillrad(g)+3\epsilon < \frac{L(g,A)}{6}$.   During the rest of the argument we will associate boundary simplices to themselves without explicitly mentioning this special case.  The arguments for these simplices will always follow from the other arguments along with the fact that the diameters are bounded by $\epsilon$.   

{\bf Step1: The 0 and 1 Skeletons}

We define the 0-chains and 1-chains, $c^0_i$ and $c^1_i$, associated to the 0-skeleton, $S_0(\Sigma)$, and the 1-skeleton, $S_1(\Sigma)$ of $\Sigma$ just as before; i.e. by mapping each 0 simplex to a closest point on $S^3$, and mapping each edge in the one skeleton to a minimizing geodesic between the endpoints.  Hence the length of the image of a 1 simplex is less than $2Fillrad(g)+3\epsilon< \frac{L(g,A)}{6}$.  We can assume (by small moves) that 0-chains of distinct vertices of $\Sigma$ are distinct and that the geodesic segments only intersect each other at endpoints.

{\bf Step 2: The 2 skeleton}

Now consider a two simplex $\Delta^2_j$ of $\Sigma$.  Let $\bar G^2_j$ be the embedded geodesic graph in $S^3$ which is the union of the geodesic segments that correspond to the edges of $G^2_j$.  We know that the support of $\partial c^2_j$ (i.e. a geodesic triangle) does not intersect $A((\bar G^2_j))$ since in fact $\bar G^2_j$ does not intersect $A((\bar G^2_j))$, for if so there would be an $x\in \bar G^2_j$ such that $A(x)$ is also in $\bar G^2_j$.  But this can't happen because the diameter of $\bar G^2_j$ is $< \frac {3} {2} \frac{L(g,A)}{6}<L(g,A).$ 

We claim that $\partial c^2_j$ represents zero in $H_1(S^3-A(\bar G^2_k))$.  Alexander duality along with the fact that $H_1(\bar G^2_j)$ is generated by $\{\partial c^2_k| \Delta^2_k$ and $\Delta^2_j $ lie in a common 4-simplex$\}$ says that we need only show that ${\cal L}(\partial c^2_j,A(\partial c^2_k))=0$ for each such $c^2_k$.  So assume that ${\cal L}(\partial c^2_j,A(\partial c^2_k))\not =0$.  Since $\Delta^2_k$ and $\Delta^2_j $ lie in a common 4-simplex they share at least one vertex and hence we let $\gamma$ be the closed simplicial curve (of combinatorial length 6) which is just $\partial \Delta^2_1$ followed by $\Delta^2_k$.  We let $\bar \gamma$ be the corresponding closed piecewise geodesic curve in $\bar G^2_j$ of length $< 6 \frac{L(g,A)}{6}=L(g,A)$ (which of course does not intersect $A(\bar \gamma)$).  By the definition of $L(g,A)$ we see that ${\cal L}(\bar \gamma,A(\bar \gamma))$ is zero as are ${\cal L}(\partial c^2_j,A(\partial c^2_j))$ and ${\cal L}(\partial c^2_k,A(\partial c^2_k))$.  On the other hand 
$$0={\cal L}(\gamma,A(\gamma))={\cal L}(\partial c^2_j,A(\partial c^2_j))+{\cal L}(\partial c^2_k,A(\partial c^2_k))+{\cal L}(\partial c^2_j,A(\partial c^2_k))+{\cal L}(\partial c^2_k,A(\partial c^2_j))$$
and hence ${\cal L}(\partial c^2_j,A(\partial c^2_k))=-{\cal L}(\partial c^2_k,A(\partial c^2_j))$.  But since $A^2=id$ and $A$ preserves orientation (since it is fixed point free) we have
$${\cal L}(\partial c^2_j,A(\partial c^2_k))={\cal L}(A(\partial c^2_j),A^2(\partial c^2_k))={\cal L}(A(\partial c^2_j),\partial c^2_k)$$ and the claim follows.

Thus we can define a 2-chain $c^2_j$ whose boundary is $\partial c^2_j$ and whose support is contained in $S^3-A(\bar G^2_j)$.

{\bf Step 3: the 3 and 4 skeleta}

Let $\Delta^3_i$ be a 3-simplex in $\Sigma$.  For each 2-face $\Delta^2_k$ the support of $c^2_k$ lies in $S^3-A(\bar G^2_k)\subset S^3-A(\bar G^3_i)$ (since $\bar G^3_i \subset \bar G^2_k$). Thus we have already defined the cycle $\partial c^3_i$ in such a way that its support is in $S^3-A(\bar G^3_i)$ and since by Alexander duality $H_3(S^3-A(\bar G^3_i))=0$ we can find a chain $c^3_i$ with boundary $\partial c^3_i$ whose support also lies in $S^3-A(\bar G^3_i)$.  Thus we can extend to the 3 skeleton.  

Now for $\Delta^4_i$ a 4-simplex of $\Sigma$ we have now defined the cycle $\partial c^4_i$ in such a way that its support is in $S^3-A(\bar G^4_i)$.  Since the support of $\partial c^4_i$ misses a point of $S^3$ there is a $c^4_i$ whose boundary is $\partial c^4_i$. 

\qed

\end{document}